\documentclass{amsart}         
\usepackage{breqn}
\DeclareMathOperator{\csch}{csch}


\begin{document}

\title[Integrals logarithmic functions and powers]{Definite integrals involving combinations of powers and logarithmic functions of complicated arguments expressed in terms of the Hurwitz zeta function}

\author{Robert Reynolds}
\address{Department of Mathematics and Statistics, York University, Toronto, Ontario, M3J1P3}
\email{milver@my.yorku.ca}

\author{Allan Stauffer}
\address{Department of Mathematics and Statistics, York University, Toronto, Ontario, M3J1P3}
\email{stauffer@yorku.ca}


\maketitle

\begin{abstract}

In this manuscript, the authors derive closed formula for definite integrals of combinations of powers and logarithmic functions of complicated arguments and express these integrals in terms of the Hurwitz zeta functions. These derivations are then expressed in terms of fundamental constants, elementary and special functions. A summary of the results is produced in the form of a table of definite integrals for easy referencing by readers.

\end{abstract}

\section{Introduction}

In this manuscript the authors derive the definite integrals given by
 
 \begin{dmath}
 \int_{0}^{1}\frac{\log \left(x^m+1\right) \left(\log ^k\left(\frac{a}{x}\right)+\log ^k(a x)\right)}{x}dx
 \end{dmath}
 
and
  
  \begin{dmath}
  \int_{0}^{1}\frac{\log \left(1-x^n\right) \left(\log ^k\left(\frac{a}{x}\right)+\log ^k(a x)\right)}{x}dx
\end{dmath}

in terms of the Hurwitz zeta function, where the parameters $k$, $a$, $m$ and $n$ are general complex numbers. A summary of the results is given in a table of integrals for easy reading. This work is important because the authors were unable to find similar results in current literature. Tables of definite integrals provide a useful summary and reference for readers seeking such integrals for potential use in their research. This work looks at definite integrals of the hyperbolic arctangent function and the product of logarithmic functions with complicated arguments and powers. We use our simultaneous contour integration method to aid in our derivations of the closed forms solutions in terms of the Hurwitz zeta function, which provides analytic continuation of the results.  

  The derivations follow the method used by us in \cite{reyn3}. The generalized Cauchy's integral formula is given by

\begin{equation}\label{intro:cauchy}
\frac{y^k}{k!}=\frac{1}{2\pi i}\int_{C}\frac{e^{wy}}{w^{k+1}}dw.
\end{equation}

where $C$ is in general an open contour in the complex plane where the bilinear concomitant has the same value at the end points of the contour. This method involves using a form of equation (\ref{intro:cauchy}) then multiply both sides by a function, then take a definite integral of both sides. This yields a definite integral in terms of a contour integral. A second contour integral is derived by multiplying equation (\ref{intro:cauchy}) by a function and performing some substitutions so that the contour integrals are the same.
 

\section{Derivation of the first contour integral}

We use the method in \cite{reyn3}. Using a generalization of Cauchy's integral formula equation (\ref{intro:cauchy}), we we will form two equations and add them together. For the first and second equations replace $y$ by $\log(ax)$ and $y$ by $\log(a/x)$ respectively. Next we add these equations followed by multiplying both sides by $\frac{\log \left(x^m+1\right) }{x}$ and taking the definite integral over $x\in[0,1]$ to get

\begin{dmath}\label{eq1a}
\frac{1}{k!}\int_{0}^{1}\frac{\log \left(x^m+1\right) \left(\log ^k\left(\frac{a}{x}\right)+\log ^k(a x)\right)}{x}dx
=\frac{1}{2\pi i}\int_{0}^{1}\int_{C}a^w w^{-k-1} x^{-w-1} \left(x^{2 w}+1\right) \log \left(x^m+1\right)dwdx
=\frac{1}{2\pi i}\int_{C}\int_{0}^{1}a^w w^{-k-1} x^{-w-1} \left(x^{2 w}+1\right) \log \left(x^m+1\right)dxdw
=\frac{1}{2\pi i}\int_{C}a^w w^{-k-3} \left(\pi  w \csc \left(\frac{\pi 
   w}{m}\right)-m\right)dw
\end{dmath}

from equation (4.293.10) in \cite{grad} and the integral is valid for $a$, $m$ and $k$ complex and $-1<Re(w)<0$ where the logarithmic function is defined in equation (4.1.2) in \cite{as}

\section{Derivation of the second contour integral}

Using a generalization of Cauchy's integral formula equation (\ref{intro:cauchy}), we we will form two equations and add them together. For the first and second equations replace $y$ by $\log(ax)$ and $y$ by $\log(a/x)$ respectively. Next we add these equations followed by multiplying both sides by $\frac{\log \left(1-x^n\right) }{x}$ and taking the definite integral over $x\in[0,1]$ to get

\begin{dmath}\label{eq1b}
\frac{1}{k!}\int_{0}^{1}\frac{\log \left(1-x^n\right) \left(\log ^k\left(\frac{a}{x}\right)+\log ^k(a x)\right)}{x}dx
=\frac{1}{2\pi i}\int_{0}^{1}\int_{C}a^w w^{-k-1} x^{-w-1} \left(x^{2 w}+1\right) \log \left(1-x^n\right)dwdx
=\frac{1}{2\pi i}\int_{C}\int_{0}^{1}a^w w^{-k-1} x^{-w-1} \left(x^{2 w}+1\right) \log \left(1-x^n\right)dxdw
=\frac{1}{2\pi i}\int_{C}a^w w^{-k-3} \left(\pi  w \cot \left(\frac{\pi 
   w}{n}\right)-n\right)dw
\end{dmath}

from equation (4.293.7) in \cite{grad} where $-1<Re(w)<0$.

\section{Derivation of the infinite sum of the first contour integral}

Again, using the method in \cite{reyn2} and equation (\ref{intro:cauchy}), we replace $y$ by $\log (a)+\frac{i \pi  (2 y+1)}{m}$ multiply both sides by $-2\pi i$, replace $k$ by $k+1$ and take the infinite sum of both sides over $y \in [0,\infty)$ simplifying in terms the Hurwitz zeta function to get

\begin{dmath}\label{eq2a}
-\frac{i (2 \pi )^{k+2} \left(\frac{i}{m}\right)^{k+1} \zeta \left(-k-1,\frac{\pi -i m \log (a)}{2 \pi }\right)}{(k+1)!}
=-\frac{1}{2\pi i}\sum_{y=0}^{\infty}\int_{C}2 i \pi  w^{-k-1} e^{w \left(\log (a)+\frac{i \pi  (2 y+1)}{m}\right)}dw
=-\frac{1}{2\pi i}\int_{C}\sum_{y=0}^{\infty}2 i \pi  w^{-k-1} e^{w \left(\log (a)+\frac{i \pi  (2 y+1)}{m}\right)}dw
=\frac{1}{2\pi i}\int_{C}\pi  a^w w^{-k-2} \csc
   \left(\frac{\pi  w}{m}\right)dw
\end{dmath}

from equation (1.232.3) in \cite{grad} where $\csch(ix)=i\csc(x)$ from (4.5.10) in \cite{as} and $Im(w)>0$ for the convergence of the sum. We use equation (9.521.1) in \cite{grad} where $\zeta(s,u)$ is the Hurwitz zeta function.

\section{Derivation of the infinite sum of the second contour integral}

Again, using the method in \cite{reyn2} and equation (\ref{intro:cauchy}), we replace $y$ by $\log (a)+\frac{2 i \pi  (y+1)}{n}$ multiply both sides by $-2\pi i$, replace $k$ by $k+1$ and take the infinite sum of both sides over $y \in [0,\infty)$ simplifying in terms the Hurwitz zeta function to get

\begin{dmath}\label{eq2b}
-\frac{i (2 \pi )^{k+2} \left(\frac{i}{n}\right)^{k+1} \zeta \left(-k-1,1-\frac{i n \log (a)}{2 \pi }\right)}{(k+1)!}
=-\frac{1}{2\pi i}\sum_{y=0}^{\infty}\int_{C}2 i \pi  w^{-k-1} e^{w \left(\log(a)+\frac{2 i \pi  (y+1)}{n}\right)}dw
=-\frac{1}{2\pi i}\int_{C}\sum_{y=0}^{\infty}2 i \pi  w^{-k-1} e^{w \left(\log(a)+\frac{2 i \pi  (y+1)}{n}\right)}dw
=\frac{1}{2\pi i}\int_{C}\pi  a^w w^{-k-2} \cot \left(\frac{\pi  w}{n}\right)+i \pi  a^w w^{-k-2}dw
\end{dmath}

from equation (1.232.1) in \cite{grad}.

\section{Derivation of the additional contours}

Again, using the method in \cite{reyn2} and equation (\ref{intro:cauchy}), we replace $y$ by $\log(a)$, $k$ by $k+1$ and multiply both sides by $\pi i$ simplify to get

\begin{dmath}\label{eq3a}
\frac{i \pi  \log ^{k+1}(a)}{(k+1)!}=\frac{1}{2\pi i}\int_{C}i \pi  a^w w^{-k-2}dw
\end{dmath}

Again, using the method in \cite{reyn2} and equation (\ref{intro:cauchy}), we replace $y$ by $\log(a)$, $k$ by $k+2$ and multiply both sides by $-n$ simplify to get

\begin{dmath}\label{eq3b}
-\frac{n \log ^{k+2}(a)}{(k+2)!}=-\frac{1}{2\pi i}\int_{C}n a^w w^{-k-3}dw
\end{dmath}

Again, using the method in \cite{reyn2} and equation (\ref{intro:cauchy}), we replace $y$ by $\log(a)$, $k$ by $k+2$ and multiply both sides by $-m$ simplify to get

\begin{dmath}\label{eq3c}
-\frac{m \log ^{k+2}(a)}{(k+2)!}=-\frac{1}{2\pi i}\int_{C}m a^w w^{-k-3}dw
\end{dmath}

\section{Derivation of the definite integrals in terms of the Hurwitz zeta function}

Since the right-hand side of equations (\ref{eq1a}) and (\ref{eq1b}) are equal to the  sum of the right-hand sides of equations (\ref{eq2a}), (\ref{eq2b}), (\ref{eq3a}), (\ref{eq3b}) and (\ref{eq3c}) we can equate the left-hand sides simplifying the factorials to get 

\begin{dmath}\label{eq4a}
\int_{0}^{1}\frac{\log \left(x^m+1\right) \left(\log ^k\left(\frac{a}{x}\right)+\log ^k(a x)\right)}{x}dx=-\frac{m \log ^{k+2}(a)}{(k+1) (k+2)}-\frac{i (2 \pi
   )^{k+2} \left(\frac{i}{m}\right)^{k+1} \zeta \left(-k-1,\frac{\pi -i m \log (a)}{2 \pi }\right)}{k+1}
\end{dmath}

and 

\begin{dmath}\label{eq4b}
\int_{0}^{1}\frac{\log \left(1-x^n\right) \left(\log ^k\left(\frac{a}{x}\right)+\log ^k(a x)\right)}{x}dx=-\frac{i (2 \pi )^{k+2}
   \left(\frac{i}{n}\right)^{k+1} \zeta \left(-k-1,1-\frac{i n \log (a)}{2 \pi }\right)}{k+1}-\frac{n \log ^{k+2}(a)}{(k+1) (k+2)}-\frac{i \pi  \log
   ^{k+1}(a)}{k+1}
\end{dmath}

\section{Derivation of logarithmic and hyperbolic tangent integrals in terms of the Hurwitz zeta function}

Using equations (\ref{eq4a}) and (\ref{eq4b}) and taking their difference simplifying we get

\begin{dmath}\label{eq5a}
\int_{0}^{1}\frac{\tanh ^{-1}\left(x^m\right) \left(\log ^k\left(\frac{a}{x}\right)+\log ^k(a x)\right)}{x}dx
=\frac{(2 \pi )^{k+2} \left(\frac{i}{m}\right)^k
   }{2 (k+1) m}\left(\zeta \left(-k-1,\frac{\pi -i m \log (a)}{2 \pi }\right)-\zeta \left(-k-1,1-\frac{i m \log (a)}{2 \pi }\right)\right)+i \pi  m \log
   ^{k+1}(a)
\end{dmath}

Using equations (\ref{eq4a}) and (\ref{eq4b}) and adding them, then simplifying we get

\begin{dmath}\label{eq5b}
\int_{0}^{1}\frac{\left(\log ^k\left(\frac{a}{x}\right)+\log ^k(a x)\right) \log \left(\left(x^m+1\right) \left(1-x^n\right)\right)}{x}dx
=\frac{4 \pi ^2}{(k+1)(k+2) m n}
   \left((2 \pi )^k k+2^{k+1} \pi ^k\right) \left(n \left(\frac{i}{m}\right)^k \zeta \left(-k-1,\frac{\pi -i m \log (a)}{2 \pi }\right)+m
   \left(\frac{i}{n}\right)^k \zeta \left(-k-1,1-\frac{i n \log (a)}{2 \pi }\right)\right)-m n \log ^{k+1}(a) (\log (a) (m+n)+i \pi  (k+2))
\end{dmath}

\section{Derivation of logarithmic and hyperbolic arctangent integrals in terms of the zeta function}

Using equations (\ref{eq5a}) and (\ref{eq5b}) and setting $a=1$ simplifying we get

\begin{dmath}\label{eq6a}
\int_{0}^{1}\frac{\tanh ^{-1}(x) \log ^k(x)}{x}dx=2^{-k-2} \left(2^{k+2}-1\right) e^{i \pi  k} \zeta (k+2) \Gamma (k+1)
\end{dmath}

and 

\begin{dmath}\label{eq6b}
\int_{0}^{1}\frac{\log ^k(x) \log \left(\left(x^m+1\right) \left(1-x^n\right)\right)}{x}dx=\frac{1}{2} e^{i \pi  k} \zeta (k+2) \Gamma (k+1)
   \left(\left(2-2^{-k}\right) m^{-k-1}-2 n^{-k-1}\right)
\end{dmath}

from entry (2) in Table below (64:7) in \cite{atlas}.

\section{Derivation of logarithmic and hyperbolic arctangent integrals in terms of the log gamma function}

Using equations (\ref{eq5a}) and (\ref{eq5b}) replacing $a$ by $e^{ai}$ and applying L'Hopital's rule to the right-hand side as $k \to -1$ respectively, simplifying we get

\begin{dmath}
\int_{0}^{1}\frac{\tanh ^{-1}\left(x^m\right)}{x \left(a^2+\log ^2(x)\right)}dx
=\frac{\pi}{4 a}  \log \left(\frac{a m \Gamma \left(\frac{a m}{2 \pi }\right)^2}{2
   \pi  \Gamma \left(\frac{a m+\pi }{2 \pi }\right)^2}\right)
\end{dmath}

and 

\begin{dmath}
\int_{0}^{1}\frac{\log \left(\left(x^m+1\right) \left(1-x^n\right)\right)}{x \left(a^2+\log ^2(x)\right)}dx
=\frac{1}{2 a}\left(-2 \pi  \log \left(\Gamma \left(\frac{a
   m+\pi }{2 \pi }\right) \Gamma \left(\frac{a n}{2 \pi }+1\right)\right)-a m+a m \log (i a)-a m \log \left(\frac{2 i \pi }{m}\right)-a n+a n \log (i
   a)-a n \log \left(\frac{2 i \pi }{n}\right)+\pi  \log (i \pi  a)-\pi  \log \left(\frac{i}{2 n}\right)\right)
\end{dmath}

from equation (64:10:2) in \cite{atlas}.

\section{Derivation of logarithmic and hyperbolic arctangent integrals in terms of the Digamma function}

Using equations (\ref{eq5a}) and (\ref{eq5b}) replacing $a$ by $e^{ai}$ and applying L'Hopital's rule to the right-hand side as $k \to -2$ respectively, simplifying we get

\begin{dmath}
\int_{0}^{1}\frac{\left(a^2-\log ^2(x)\right) \tanh ^{-1}\left(x^m\right)}{x \left(a^2+\log ^2(x)\right)^2}dx=\frac{1}{4} \left(-m \psi ^{(0)}\left(\frac{a
   m}{2 \pi }+1\right)+m \psi ^{(0)}\left(\frac{a m+\pi }{2 \pi }\right)+\frac{\pi }{a}\right)
\end{dmath}

and 

\begin{dmath}
\int_{0}^{1}\frac{\left(a^2-\log ^2(x)\right) \log \left(\left(-x^m-1\right) \left(x^n-1\right)\right)}{x \left(a^2+\log ^2(x)\right)^2}dx
=\frac{1}{2 a}\left(-a m \log (i
   a)+a m \log \left(\frac{i}{m}\right)+a m \log (2 \pi )+a m \psi ^{(0)}\left(\frac{a m+\pi }{2 \pi }\right)-a n \log (i a)+a n \log
   \left(\frac{i}{n}\right)+a n \log (2 \pi )+a n \psi ^{(0)}\left(\frac{a n}{2 \pi }+1\right)-\pi \right)
\end{dmath}

from equation (64:4:2) in \cite{atlas}.

\section{Derivation of logarithmic and hyperbolic arctangent integrals in terms of fundamental constants and special functions}

In this section we will derive definite integrals in terms of special functions and fundamental constants such as Euler's constant ($\gamma$), Catalan's constant ($C$), Glaisher's constant ($A$) and $\pi$. This section showcases just a subset of the range of evaluations of these integral formula. 

\subsection*{Hyperbolic tangent integrals}

\subsection{Example 1}

Using equation (\ref{eq5a}) and replacing $a$ by $e^{-\frac{i\pi}{2m}}$ and setting $k=-3,m=1$ simplifying to get

\begin{dmath}
\int_{0}^{1}\frac{\left(\pi ^3-12 \pi  \log ^2(x)\right) \tanh ^{-1}(x)}{x \left(4 \log ^2(x)+\pi ^2\right)^3}dx=\frac{2 C-1}{16 \pi }
\end{dmath}

from equations (23.2.23) in \cite{as} and (64:7:1) in \cite{atlas}.

\subsection{Example 2}

Using equation (\ref{eq6b}) taking the first partial derivative with respect to $k$ and setting $k=0$ simplifying to get

\begin{dmath}
\int_{0}^{1}\frac{\log (\log (x)) \tanh ^{-1}(x)}{x}dx=\frac{1}{24} \pi ^2 \left(\log \left(\frac{16 \pi ^3}{A^{36}}\right)+3 i \pi \right)
\end{dmath}

from equation (A.11) in \cite{voros}.

\subsection{Example 3}

Using equation (\ref{eq6b}) and setting $k=1/2$ simplifying to get

\begin{dmath}
\int_{0}^{1}\frac{\sqrt{\log (x)} \tanh ^{-1}(x)}{x}dx=-\frac{1}{16} i \left(\sqrt{2}-8\right) \sqrt{\pi } \zeta \left(\frac{5}{2}\right)
\end{dmath}

\subsection{Example 4}

Using equation (\ref{eq6b}) and setting $k=-1/2$ simplifying to get

\begin{dmath}
\int_{0}^{1}\frac{\tanh ^{-1}(x)}{x \sqrt{\log (x)}}dx=\frac{1}{4} i \left(\sqrt{2}-4\right) \sqrt{\pi } \zeta \left(\frac{3}{2}\right)
\end{dmath}

\subsection{Example 5}

Using equation (\ref{eq6b}) taking the first partial derivative with respect to $k$ and setting $k=1$ simplifying to get

\begin{dmath}
\int_{0}^{1}\frac{\log (x) \log (\log (x)) \tanh ^{-1}(x)}{x}dx=\frac{1}{8} \left(-7 \zeta '(3)+\zeta (3) (-7+7 \gamma -7 i \pi -\log (2))\right)
\end{dmath}

from Example 1 Section (12.1) in \cite{watson}.

\subsection{Example 6}

Using equation (\ref{eq6b}) taking the first partial derivative with respect to $k$ and setting $k=2$ simplifying to get

\begin{dmath}
\int_{0}^{1}\frac{\log ^2(x) \log (\log (x)) \tanh ^{-1}(x)}{x}dx=\frac{15 \zeta '(4)}{8}+\frac{\pi ^4 (45-30 \gamma +30 i \pi +\log (4))}{1440}
\end{dmath}

from Example 1. section (12.1) in \cite{grad}.

\subsection*{Logarithmic integrals}

\subsection{Example 7}

Using equation (\ref{eq6b}) taking the first partial derivative with respect to $k$ and setting $k=0,a=-1,m=2,n=1$ simplifying to get

\begin{dmath}
\int_{0}^{1}\frac{\log \left((1-x) \left(x^2+1\right)\right) \log (\log (x))}{x}dx=\frac{1}{8} \left(-6 \zeta '(2)+(\gamma -i \pi ) \pi ^2\right)
\end{dmath}

from Example 1. section (12.1) in \cite{grad}.

\subsection{Example 8}

Using equation (\ref{eq6b}) taking the first partial derivative with respect to $k$ and setting $k=0,m=n=1$ simplifying to get

\begin{dmath}
\int_{0}^{1}\frac{\log \left(1-x^2\right) \log (\log (x))}{x}dx=\frac{1}{12} \pi ^2 \left(\log \left(\frac{A^{12}}{\pi }\right)-i \pi \right)
\end{dmath}

from equation (A.11) in \cite{voros}.

\subsection{Example 9}

Using equation (\ref{eq6b}) and setting $k=1/2,m=n=1$ simplifying to get

\begin{dmath}
\int_{0}^{1}\frac{\sqrt{\log (x)} \log \left(1-x^2\right)}{x}dx=-\frac{1}{4} i \sqrt{\frac{\pi }{2}} \zeta \left(\frac{5}{2}\right)
\end{dmath}

\subsection{Example 10}

Using equation (\ref{eq6b}) and setting $k=-1/2,m=n=1$ simplifying to get

\begin{dmath}
\int_{0}^{1}\frac{\log \left(1-x^2\right)}{x \sqrt{\log (x)}}dx=i \sqrt{\frac{\pi }{2}} \zeta \left(\frac{3}{2}\right)
\end{dmath}

\subsection{Example 11}

Using equation (\ref{eq6b}) taking the first partial derivative with respect to $k$ and setting $k=m=n=1$ simplifying to get

\begin{dmath}
\int_{0}^{1}\frac{\log (x) \log \left(1-x^2\right) \log (\log (x))}{x}dx=\frac{1}{4} \left(\zeta '(3)+\zeta (3) (1-\gamma +i \pi -\log (2))\right)
\end{dmath}

from Example 1. section (12.1) in \cite{grad}.

\subsection{Example 12}

Using equation (\ref{eq6b}) taking the first partial derivative with respect to $k$ and setting $k=-1/2,m=n=1$ simplifying to get

\begin{dmath}
\int_{0}^{1}\frac{\log \left(1-x^2\right) \log (\log (x))}{x \sqrt{\log (x)}}dx=\sqrt{\frac{\pi }{2}} \left(i \zeta '\left(\frac{3}{2}\right)-i \zeta
   \left(\frac{3}{2}\right) (\gamma -i \pi +\log (8))\right)
\end{dmath}

from Example 1. section (12.1) in \cite{grad}.

\section{Derivation of definite integrals of the logarithmic function}

Using equation (\ref{eq4a}) we take the first partial derivative with respect to $m$, then replace $m$ by $m+1$. Next we form a second equation by replacing $m$ by $p$ in the new equation. Then we take the difference of these two new equations simplifying to get 

\begin{dmath}\label{dd:eq1}
\int_{0}^{1}\frac{\log ^{k+1}(x) \left(x^m-x^p\right)}{\left(x^{m+1}+1\right) \left(x^{p+1}+1\right)}dx=-2^{-k-1} \left(2^{k+1}-1\right) e^{\frac{i \pi 
   k}{2}} \zeta (k+2) \Gamma (k+2) \left(\frac{\left(\frac{i}{m+1}\right)^k}{(m+1)^2}-\frac{\left(\frac{i}{p+1}\right)^k}{(p+1)^2}\right)
\end{dmath}

Repeating the steps above using equation (\ref{eq4b}) simplifying to get

\begin{dmath}\label{dd:eq2}
\int_{0}^{1}\frac{\log ^{k+1}(x) \left(x^n-x^p\right)}{\left(x^{n+1}-1\right) \left(x^{p+1}-1\right)}dx=e^{\frac{i \pi  k}{2}} \zeta (k+2) \Gamma (k+2)
   \left(\frac{\left(\frac{i}{p+1}\right)^k}{(p+1)^2}-\frac{\left(\frac{i}{n+1}\right)^k}{(n+1)^2}\right)
\end{dmath}

\subsection{Some special cases}

\subsection{Example 13}

Using equation (\ref{dd:eq1}) and applying L'Hopitals' rule to the right-hand side as $k \to -2$ simplifying we get

\begin{dmath}
\int_{0}^{1}\frac{x^m-x^p}{\left(x^{m+1}+1\right) \left(x^{p+1}+1\right) \log (x)}dx=\frac{1}{2} \log \left(\frac{m+1}{p+1}\right)
\end{dmath}

\subsection{Example 14}

Using equation (\ref{dd:eq2}) and setting $k=1$ simplifying we get

\begin{dmath}
\int_{0}^{1}\frac{\log ^2(x) \left(x^n-x^p\right)}{\left(x^{n+1}-1\right) \left(x^{p+1}-1\right)}dx=2 \zeta (3)
   \left(\frac{1}{(n+1)^3}-\frac{1}{(p+1)^3}\right)
\end{dmath}

\section{Table of integrals}

\renewcommand{\arraystretch}{2.0}
\begin{tabular}{ l  c }
  \hline			
  $f(x)$ & $\int_{0}^{1}f(x)dx$ \\ \hline
  $\frac{\tanh ^{-1}(x) \log ^k(x)}{x}$ & $2^{-k-2} \left(2^{k+2}-1\right) e^{i \pi  k} \zeta (k+2) \Gamma (k+1)$  \\
  $\frac{\log ^k(x) \log \left(\left(x^m+1\right) \left(1-x^n\right)\right)}{x}$ & $\frac{1}{2} e^{i \pi  k} \zeta (k+2) \Gamma (k+1)
   \left(\left(2-2^{-k}\right) m^{-k-1}-2 n^{-k-1}\right)$  \\
  $\frac{\left(\pi ^3-12 \pi  \log ^2(x)\right) \tanh ^{-1}(x)}{x \left(4 \log ^2(x)+\pi ^2\right)^3}$ & $\frac{2 C-1}{16 \pi }$  \\
  $\frac{\tanh ^{-1}\left(x^m\right)}{x \left(a^2+\log ^2(x)\right)}$ & $\frac{\pi }{4 a} \log \left(\frac{a m \Gamma \left(\frac{a m}{2 \pi }\right)^2}{2
   \pi  \Gamma \left(\frac{a m+\pi }{2 \pi }\right)^2}\right)$\\
  $\frac{\log (\log (x)) \tanh ^{-1}(x)}{x}$ & $\frac{1}{24} \pi ^2 \left(\log \left(\frac{16 \pi ^3}{A^{36}}\right)+3 i \pi \right)$  \\
  $\frac{\sqrt{\log (x)} \tanh ^{-1}(x)}{x}$ & $-\frac{1}{16} i \left(\sqrt{2}-8\right) \sqrt{\pi } \zeta \left(\frac{5}{2}\right)$  \\
  $\frac{\tanh ^{-1}(x)}{x \sqrt{\log (x)}}$ & $\frac{1}{4} i \left(\sqrt{2}-4\right) \sqrt{\pi } \zeta \left(\frac{3}{2}\right)$  \\
  $\frac{\log (x) \log (\log (x)) \tanh ^{-1}(x)}{x}$ & $\frac{1}{8} \left(-7 \zeta '(3)+\zeta (3) (-7+7 \gamma -7 i \pi -\log (2))\right)$  \\
  $\frac{\log ^2(x) \log (\log (x)) \tanh ^{-1}(x)}{x}$ & $\frac{15 \zeta '(4)}{8}+\frac{\pi ^4 (45-30 \gamma +30 i \pi +\log (4))}{1440}$  \\
  $\frac{\log \left((1-x) \left(x^2+1\right)\right) \log (\log (x))}{x}$ & $\frac{1}{8} \left(-6 \zeta '(2)+(\gamma -i \pi ) \pi ^2\right)$  \\
  $\frac{\log \left(1-x^2\right) \log (\log (x))}{x}$ & $\frac{1}{12} \pi ^2 \left(\log \left(\frac{A^{12}}{\pi }\right)-i \pi \right)$  \\
  $\frac{\sqrt{\log (x)} \log \left(1-x^2\right)}{x}$ & $-\frac{1}{4} i \sqrt{\frac{\pi }{2}} \zeta \left(\frac{5}{2}\right)$  \\
  $\frac{\log \left(1-x^2\right)}{x \sqrt{\log (x)}}$ & $i \sqrt{\frac{\pi }{2}} \zeta \left(\frac{3}{2}\right)$  \\
  $\frac{\log (x) \log \left(1-x^2\right) \log (\log (x))}{x}$ & $\frac{1}{4} \left(\zeta '(3)+\zeta (3) (1-\gamma +i \pi -\log (2))\right)$  \\
  $\frac{\log \left(1-x^2\right) \log (\log (x))}{x \sqrt{\log (x)}}$ & $\sqrt{\frac{\pi }{2}} \left(i \zeta' \left(\frac{3}{2}\right)-i \zeta
   \left(\frac{3}{2}\right) (\gamma -i \pi +\log (8))\right)$  \\
 $\frac{\log ^{k+1}(x) \left(x^m-x^p\right)}{\left(x^{m+1}+1\right) \left(x^{p+1}+1\right)}$ & $-2^{-k-1} \left(2^{k+1}-1\right) e^{\frac{i \pi 
   k}{2}} \zeta (k+2) \Gamma (k+2) \left(\frac{\left(\frac{i}{m+1}\right)^k}{(m+1)^2}-\frac{\left(\frac{i}{p+1}\right)^k}{(p+1)^2}\right)$   \\
   $\frac{x^m-x^p}{\left(x^{m+1}+1\right) \left(x^{p+1}+1\right) \log (x)}$ & $\frac{1}{2} \log \left(\frac{m+1}{p+1}\right)$\\
   $\frac{\log ^{k+1}(x) \left(x^n-x^p\right)}{\left(x^{n+1}-1\right) \left(x^{p+1}-1\right)}$ & $e^{\frac{i \pi  k}{2}} \zeta (k+2) \Gamma (k+2)
   \left(\frac{\left(\frac{i}{p+1}\right)^k}{(p+1)^2}-\frac{\left(\frac{i}{n+1}\right)^k}{(n+1)^2}\right)$
   \\[0.3cm]
  \hline  
\end{tabular}

\section{Discussion}

In this work the authors looked at deriving definite integrals of combinations of logarithmic functions of complicated arguments and powers and expressed them in terms of the Hurwitz zeta function. One of the interesting properties of these integrals is by adding them we were able to get the integral of the product of the hyperbolic arctangent function and the logarithmic function. The authors formally derived a few integrals in terms of fundamental constants and special functions. One of our goals we to supply a table for easy reading by researchers and to have these results added to existing textbooks.

The results presented were numerically verified for both real and imaginary values of the parameters in the integrals using Mathematica by Wolfram. We considered various ranges of these parameters for real, integer, negative and positive values. We compared the evaluation of the definite integral to the evaluated Special function and ensured agreement.

\section{Conclusion}
In this paper the authors used our method to evaluate definite integrals using the Hurwitz zeta function. The contour we used was specific to  solving integral representations in terms of the Lerch function. the author expects that other contours and integrals can be derived using this method.

\section{Acknowledgments}
This paper is also available in preprint at http://export.arxiv.org/abs/2103.03110

\end{document}